\documentclass{elsarticle}
\usepackage{amsmath}
\usepackage{amssymb}
\usepackage{gastex}
\usepackage{graphicx}
\usepackage{kbordermatrix}

\newtheorem{Theorem}{Theorem}
\newtheorem{Lemma}[Theorem]{Lemma}
\newtheorem{Proposition}[Theorem]{Proposition}
\newtheorem{Corollary}[Theorem]{Corollary}
\newdefinition{Definition}[Theorem]{Definition}
\newdefinition{Example}[Theorem]{Example}
\newdefinition{Remark}[Theorem]{Remark}
\newproof{proof}{Proof}
\newenvironment{Proof}{\begin{proof}}{\qed\end{proof}}

\newcommand{\kernel}{\ker}
\newcommand{\pair}[2]{\{{#1},{#2}\}}
\newcommand{\two}{\mathbb{F}_2}
\newcommand{\sub}[1]{[ #1 ]}
\newcommand{\xor}{\oplus}
\newcommand{\gnrdom}{\mathrm{gnrdom}}
\newcommand{\snrdom}{\mathrm{snrdom}}

\renewcommand{\emptyset}{\varnothing}
\renewcommand{\kbldelim}{(}
\renewcommand{\kbrdelim}{)}

\begin{document}

%
%

\begin{frontmatter}

\title{Maximal Pivots on Graphs\\ with an Application to Gene Assembly}

\author{Robert Brijder\corref{cor}}
\cortext[cor]{Corresponding author} \ead{rbrijder@liacs.nl}

\author{Hendrik Jan Hoogeboom}

\address{Leiden Institute of Advanced Computer Science,\\
Leiden University, The Netherlands}

\begin{abstract}
We consider principal pivot transform (\emph{pivot}) on graphs. We
define a natural variant of this operation, called dual pivot, and
show that both the kernel and the set of maximally applicable pivots
of a graph are invariant under this operation. The result is
motivated by and applicable to the theory of gene assembly in
ciliates.
\end{abstract}

\begin{keyword}
principal pivot transform \sep algebraic graph theory \sep overlap
graph \sep gene assembly in ciliates
\end{keyword}

\end{frontmatter}

\section{Introduction}
The pivot operation, due to Tucker \cite{tucker1960}, partially
(component-wise) inverts a given matrix. It appears naturally in
many areas including mathematical programming and numerical
analysis, see \cite{Tsatsomeros2000151} for a survey. Over $\two$
(which is the natural setting to consider for graphs), the pivot
operation has, in addition to matrix and graph interpretations
\cite{Geelen97}, also an interpretation in terms of delta matroids
\cite{bouchet1987}.

In this paper we define the \emph{dual pivot}, which has an
identical effect on graphs as the (regular) pivot, however the
condition for it to be applicable differs. The main result of the
paper is that any two graphs in the same orbit under dual pivot have
the same family of maximal pivots (cf. Theorem~\ref{thm:dual}),
i.e., the same family of maximally partial inverses of that matrix.
This result is obtained by combining each of the aforementioned
interpretations of pivot.

This research is motivated by the theory of gene assembly in
ciliates \cite{GeneAssemblyBook}, which is recalled in
Section~\ref{sec:appl_ga}. Without the context of gene assembly this
main result (Theorem~\ref{thm:dual}) is surprising; it is not found
in the extensive literature on pivots. It fits however with the
intuition and results from the string based model of gene assembly
\cite{FibersRangeRedGraphs/Brijder07}, and in this paper we
formulate it for the more general graph based model. It is
understood and proven here using completely different techniques,
algebraical rather than combinatorial.

\section{Notation and Terminology} \label{sec:notation}
The field with two elements is denoted by $\two$. Our matrix
computations will be over $\two$. Hence addition is equal to the
logical exclusive-or, also denoted by $\xor$, and multiplication is
equal to the logical conjunction, also denoted by $\land$. These
operations carry over to sets, e.g., for sets $A, B \subseteq V$ and
$x \in V$, $x \in A \xor B$ iff $(x \in A) \xor (x \in B)$.

A set system is a tuple $M=(V,D)$, where $V$ is a finite set and $D
\subseteq 2^V$ is a set of subsets of $V$. Let $\min(D)$ ($\max(D)$,
resp.) be the family of minimal (maximal, resp.) sets in $D$ w.r.t.
set inclusion, and let $\min(M) = (V,\min(D))$ ($\max(M) =
(V,\max(D))$, resp.) be the corresponding set systems.

Let $V$ be a finite set, and $A$ be a $V \times V$-matrix (over an
arbitrary field), i.e., $A$ is a matrix where the rows and columns
of $A$ are identified by elements of $V$. Therefore, e.g., the
following matrices with $V = \{p,q\}$ are equal: $ \kbordermatrix{
  & p & q \\
p & 1 & 1 \\
q & 0 & 1}$ and $\kbordermatrix{
  & q & p \\
q & 1 & 0 \\
p & 1 & 1}$. For $X \subseteq V$, the principal submatrix of $A$
w.r.t. $X$ is denoted by $A[X]$, i.e., $A[X]$ is the $X \times
X$-matrix obtained from $A$ by restricting to rows and columns in
$X$. Similarly, we define $A \backslash X = A[V \backslash X]$.
Notions such as matrix inversion $A^{-1}$ and determinant $\det(A)$
are well defined for $V \times V$-matrices. By convention, $\det
(A\sub{\varnothing}) = 1$.

A set $X \subseteq V$ is called dependent in $A$ iff the columns of
$A$ corresponding to $X$ are linearly dependent. We define
$\mathcal{P}_A = (I,D)$ to be the partition of $2^{V}$ such that $D$
($I$, respectively) contains the dependent (independent,
respectively) subsets of $V$ in $A$. By convention, $\emptyset \in
I$. The sets in $\max(I)$ are called the bases of $A$.

We have that $\mathcal{P}_A = (I,D)$ is uniquely determined by
$\max(I)$ (and the set $V$). Similarly, $\mathcal{P}_A$ is uniquely
determined by $\min(D)$ (and the set $V$). These properties are
specifically used in matroid theory, where a matroid may be
described by its independent sets $(V,I)$, by its family of bases
$(V,\max(I))$, or by its circuits $(V,\min(D))$. Moreover, for each
basis $X \in \max(I)$, $|X|$ is equal to the rank $r$ of $A$.

We consider undirected graphs without parallel edges, however we do
allow loops.
For a graph $G=(V,E)$ we use $V(G)$ and $E(G)$ to denote its set of
vertices $V$ and set of edges $E$, respectively, where for $x \in
V$, $\{x\} \in E$ iff $x$ has a loop. For $X \subseteq V$, we denote
the subgraph of $G$ induced by $X$ as $G\sub{X}$.

With a graph $G$ one associates its adjacency matrix $A(G)$, which
is a $V \times V$-matrix $\left(a_{u,v}\right)$ over $\two$ with
$a_{u,v} = 1$ iff $\{u,v\} \in E$.
The matrices corresponding to graphs are precisely the symmetric
$\two$-matrices; loops corresponding to diagonal $1$'s. Note that
for $X \subseteq V$, $A(G\sub{X}) = (A(G))[X]$.

Over $\two$, vectors indexed by $V$ can be identified with subsets
of $V$, and a $V \times V$-matrix defines a linear transformation on
subsets of $V$.
The kernel (also called null space) of a matrix $A$, denoted by
$\kernel(A)$ is determined by those linear combinations of column
vectors of $A$ that sum up to the zero vector $0$. Working in
$\two$, we regard the elements of $\kernel(A)$ as subsets of $V$.
Moreover, the kernel of $A$ is the eigenspace $E_0(A)$ on value $0$,
and similar as $\kernel(A)$, the elements of the (only other)
eigenspace $E_1(A) = \{v \in V \mid Av = v \}$ on value $1$ are also
considered as sets.

We will often identify a graph with its adjacency matrix, so, e.g.,
by the determinant of graph $G$, denoted by $\det G$, we will mean
the determinant $\det A(G)$ of its adjacency matrix computed over
$\two$. In the same vein we will often simply write $\kernel(G)$,
$E_1(G)$, $\mathcal{P}_G$, etc.

Let $\mathcal{P}_G = (I,D)$ for some graph $G$. As $G$ is a $V
\times V$-matrix over $\two$, we have that $X \in D$ iff there is a
$S \subseteq X$ with $S \in \kernel(G) \backslash \{\varnothing\}$.
Moreover, $\min(D) = \min(\kernel(G)\backslash \{\varnothing\})$ and
$\kernel(G)$ is the closure of $\min(D)$ under $\xor$ (i.e.,
$\min(D)$ spans $\kernel(G)$). Consequently, $\min(D)$ uniquely
determines $\kernel(G)$ and vice versa. As $\min(D)$ in turn
uniquely determines $\mathcal{P}_G$, the following holds.

\begin{Corollary} \label{cor:ker_iff_max_indep}
For graphs $G_1$ and $G_2$, $\kernel(G_1) = \kernel(G_2)$ iff the
families of bases of $G_1$ and of $G_2$ are equal.
\end{Corollary}

\section{Pivots} \label{sec:def_pivots}

In general the pivot operation can be studied for matrices over
arbitrary fields, e.g., as done in \cite{Tsatsomeros2000151}. In
this paper we restrict ourselves to symmetric matrices over $\two$,
which leads to a number of additional viewpoints to the same
operation, and for each of them an equivalent definition for
pivoting. Each of these definitions is known, but (to our best
knowledge) they were not before collected in one text.

\paragraph{Matrices}
Let $A$ be a $V \times V$-matrix (over an arbitrary field), and let
$X \subseteq V$ be such that $A\sub{X}$ is nonsingular, i.e., $\det
A\sub{X} \neq 0$. The \emph{pivot} of $A$ on $X$, denoted by $A*X$,
is defined as follows, see \cite{tucker1960}. Let $A = \left(
\begin{array}{c|c}
P & Q \\
\hline R & S
\end{array}
\right)$ with $P = A\sub{X}$. Then
$$
A*X = \left(
\begin{array}{c|c}
P^{-1} & -P^{-1} Q \\
\hline R P^{-1} & S - R P^{-1} Q
\end{array}
\right).
$$
Matrix $(A*X)\setminus X = S - R P^{-1} Q$ is called the \emph{Schur
complement} of $X$ in $A$.

The pivot is sometimes considered a partial inverse, as $A$ and
$A*X$ are related by the following characteristic equality, where
the vectors $x_1$ and $y_1$ correspond to the elements of $X$. In
fact, this formula defines $A*X$ given $A$ and $X$
\cite{Tsatsomeros2000151}.
\begin{eqnarray} \label{pivot_def_reverse} A \left(
\begin{array}{c} x_1 \\ x_2 \end{array} \right) = \left(\begin{array}{c} y_1 \\ y_2 \end{array} \right) \mbox{ iff } A*X \left(
\begin{array}{c} y_1 \\ x_2 \end{array} \right) = \left(\begin{array}{c} x_1 \\ y_2 \end{array}
\right) \end{eqnarray} Note that if $\det A \not= 0$, then $A * V =
A^{-1}$. By Equation~(\ref{pivot_def_reverse}) we see that a pivot
operation is an involution (i.e., operation of order $2$), and more
generally, if $(A*X)*Y$ is defined, then $A*(X \xor Y)$ is defined
and they are equal.

The following fundamental result on pivots is due to
Tucker~\cite{tucker1960} (see also
\cite[Theorem~4.1.1]{cottle1992}). It is used in
\cite{BHH/PivotsDetPM/09} to study sequences of pivots.

\begin{Proposition}[\cite{tucker1960}]\label{prop:tucker}
Let $A$ be a $V \times V$-matrix, and let $X\subseteq V$ be such
that $\det A\sub{X} \neq 0$. Then, for $Y \subseteq V$, $\det
(A*X)\sub{Y} = \det A\sub{X \xor Y} / \det A\sub{X}$.
\end{Proposition}
It may be interesting to remark here that
Proposition~\ref{prop:tucker} for the case $Y = V \setminus X$ is
called the Schur determinant formula and was shown already in 1917
by Issai Schur, see \cite{Schur1917detformula}.

It is easy to verify from the definition of pivot that $A*X$ is
skew-symmetric whenever $A$ is. In particular, if $G$ is a graph
(i.e., a symmetric matrix over $\two$), then $G*X$ is also a graph.
From now on we restrict our attention to graphs.

\paragraph{Delta Matroids}
Consider now a set system $M = (V,D)$. We define, for $X \subseteq
V$, the \emph{twist} $M * X = (V,D * X)$, where $D * X = \{Y \xor X
\mid Y \in D\}$.

Let $G$ be a graph and let $\mathcal{M}_G = (V(G),D_G)$ be the set
system with $D_G = \{ X \subseteq V(G) \mid \det G\sub{X} = 1\}$. It
is easy to verify that $G$ can be (re)constructed given
$\mathcal{M}_G$: $\{u\}$ is a loop in $G$ iff $\{u\} \in D_G$, and
$\{u,v\}$ is an edge in $G$ iff $(\{u,v\} \in D_G) \xor ((\{u\} \in
D_G) \wedge (\{v\} \in D_G))$, see
\cite[Property~3.1]{Bouchet_1991_67}. In this way, the family of
graphs (with set $V$ of vertices) can be considered as a subset of
the family of set systems (over set $V$).

Proposition~\ref{prop:tucker} allows for another (equivalent)
definition of pivot over $\two$. Indeed, over $\two$, we have by
Proposition~\ref{prop:tucker}, $\det (A*X)\sub{Y} = \det A\sub{X
\xor Y}$ for all $Y \subseteq V$ assuming $A*X$ is defined.
Therefore, for $\mathcal{M}_{G*X}$ we have $D_{G*X} = \{ Z \mid \det
((G*X)\sub{Z}) = 1\} = \{ Z \mid \det (G\sub{X \xor Z}) = 1\} = \{ X
\xor Y \mid \det (G\sub{Y}) = 1\} = D_G * X$, see
\cite{bouchet1987}. Hence $\mathcal{M}_G
* X = \mathcal{M}_{G*X}$ is an alternative definition of the pivot
operation over $\two$.

It turns out that $\mathcal{M}_G$ has a special structure, that of a
\emph{delta matroid}, allowing a specific exchange of elements
between any two sets of $D_G$, see \cite{bouchet1987}. However, not
every delta matroid $M$ has a graph representation, i.e., $M$ may
not be of the form $\mathcal{M}_G$ for any graph $G$ (a
characterization of such representable delta matroids over $\two$ is
given in \cite{Bouchet_1991_67}).

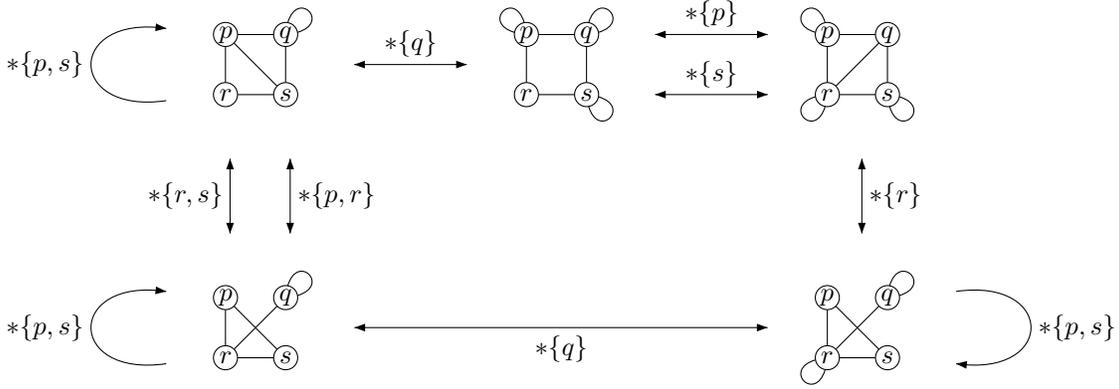
\begin{figure}[t]
\unitlength 1mm%
\begin{center}
{\begin{picture}(100,45)(-10,-05)
\node[Nw=25,Nh=25,Nframe=n](I)(00,35){%
  \unitlength0.4mm
  \begin{picture}(30,30)
  \gasset{AHnb=0,Nw=1.5,Nh=1.5,Nframe=n,Nfill=y}
  \gasset{AHnb=0,Nw=8,Nh=8,Nframe=y,Nfill=n}
    \node(r)(05,05){$r$}
    \node(p)(05,25){$p$}
    \node(q)(25,25){$q$}
    \node(s)(25,05){$s$}
    \drawedge(p,q){}
    \drawedge(p,r){}
    \drawedge(p,s){}
    \drawedge(q,s){}
    \drawedge(r,s){}
    \drawloop[loopangle=45,loopdiam=7](q){}
  \end{picture}
}
\node[Nw=25,Nh=25,Nframe=n](II)(40,35){%
  \unitlength0.4mm
  \begin{picture}(30,30)
  \gasset{AHnb=0,Nw=1.5,Nh=1.5,Nframe=n,Nfill=y}
  \gasset{AHnb=0,Nw=8,Nh=8,Nframe=y,Nfill=n}
  \node(r)(05,05){$r$}
  \node(p)(05,25){$p$}
  \node(q)(25,25){$q$}
  \node(s)(25,05){$s$}
  \drawedge(p,q){}
  \drawedge(p,r){}
  \drawedge(q,s){}
  \drawedge(r,s){}
  \drawloop[loopangle=45,loopdiam=7](q){}
  \drawloop[loopangle=135,loopdiam=7](p){}
  \drawloop[loopangle=-45,loopdiam=7](s){}
\end{picture}
}
\node[Nw=25,Nh=25,Nframe=n](III)(80,35){%
  \unitlength0.4mm
  \begin{picture}(30,30)
  \gasset{AHnb=0,Nw=1.5,Nh=1.5,Nframe=n,Nfill=y}
  \gasset{AHnb=0,Nw=8,Nh=8,Nframe=y,Nfill=n}
  \node(r)(05,05){$r$}
  \node(p)(05,25){$p$}
  \node(q)(25,25){$q$}
  \node(s)(25,05){$s$}
  \drawedge(p,q){}
  \drawedge(p,r){}
  \drawedge(q,s){}
  \drawedge(r,s){}
  \drawedge(q,r){}
  \drawloop[loopangle=225,loopdiam=7](r){}
  \drawloop[loopangle=135,loopdiam=7](p){}
  \drawloop[loopangle=-45,loopdiam=7](s){}
\end{picture}
}
\node[Nw=25,Nh=25,Nframe=n](IV)(80,00){%
  \unitlength0.4mm
  \begin{picture}(30,30)
  \gasset{AHnb=0,Nw=1.5,Nh=1.5,Nframe=n,Nfill=y}
  \gasset{AHnb=0,Nw=8,Nh=8,Nframe=y,Nfill=n}
  \node(r)(05,05){$r$}
  \node(p)(05,25){$p$}
  \node(q)(25,25){$q$}
  \node(s)(25,05){$s$}
  \drawedge(p,r){}
  \drawedge(p,s){}
  \drawedge(q,r){}
  \drawedge(r,s){}
  \drawloop[loopangle=45,loopdiam=7](q){}
  \drawloop[loopangle=225,loopdiam=7](r){}
\end{picture}
}
\node[Nw=25,Nh=25,Nframe=n](V)(00,00){%
  \unitlength0.4mm
  \begin{picture}(30,30)
  \gasset{AHnb=0,Nw=1.5,Nh=1.5,Nframe=n,Nfill=y}
  \gasset{AHnb=0,Nw=8,Nh=8,Nframe=y,Nfill=n}
  \node(r)(05,05){$r$}
  \node(p)(05,25){$p$}
  \node(q)(25,25){$q$}
  \node(s)(25,05){$s$}
  \drawedge(p,r){}
  \drawedge(p,s){}
  \drawedge(q,r){}
  \drawedge(r,s){}
  \drawloop[loopangle=45,loopdiam=7](q){}
\end{picture}
}
  \drawedge[AHnb=1,ATnb=1](I,II){$*\{q\}$}
  \drawloop[loopangle=180,loopdiam=10](I){$*\{p,s\}$}
  \drawedge[AHnb=1,ATnb=1,syo=4,eyo=4](II,III){$*\{p\}$}
  \drawedge[AHnb=1,ATnb=1,syo=-4,eyo=-4](II,III){$*\{s\}$}
  \drawedge[AHnb=1,ATnb=1](III,IV){$*\{r\}$}
  \drawloop[loopangle=00,loopdiam=10](IV){$*\{p,s\}$}
  \drawedge[AHnb=1,ATnb=1](IV,V){$*\{q\}$}
  \drawedge[AHnb=1,ATnb=1,sxo=4,exo=4,ELside=r](V,I){$*\{p,r\}$}
  \drawedge[AHnb=1,ATnb=1,sxo=-4,exo=-4](V,I){$*\{r,s\}$}
  \drawloop[loopangle=180,loopdiam=10](V){$*\{p,s\}$}
\end{picture}}
\end{center}
\caption{The orbit of $G$ under pivot. Only the elementary pivots
are shown.}\label{fig:pivot_space}
\end{figure}

\begin{Example}\label{ex:introd_pivot}
Let $G$ be the graph depicted in the upper-left corner of
Figure~\ref{fig:pivot_space}. We have $A(G) = \kbordermatrix{
  & p & q & r & s \\
p & 0 & 1 & 1 & 1 \\
q & 1 & 1 & 0 & 1 \\
r & 1 & 0 & 0 & 1 \\
s & 1 & 1 & 1 & 0 }$. This corresponds to $\mathcal{M}_G =
(\{p,q,r,s\},D_G)$, where
$$D_G = \{\varnothing, \{q\},
\{p,q\}, \{p,r\}, \{p,s\}, \{q,s\}, \{r,s\}, \{p,q,s\}, \{p,q,r\},
\{q,r,s\} \}.
$$
For example, $\{p,q\} \in D_G$ since $\det(G[\{p,q\}]) = \det \left(
\begin{array}{cc}
0 & 1 \\
1 & 1 \end{array} \right) = 1$. Then $D_G * \{p,q\} =
\{ \varnothing, \{p\}, \{r\}, \{s\}, \{p,q\}, \{p,s\}, \{q,r\},
\{q,s\}, \{p,r,s\}, \{p,q,r,s\} \}$,
%
%
and the corresponding graph is depicted on the top-right in the same
Figure~\ref{fig:pivot_space}. Equivalently, this graph is obtained
from $G$ by pivot on $\{p,q\}$. Also note that we have $D_G *
\{p,s\} = D_G$, and therefore the pivot of $G$ on $\{p,s\}$ obtains
$G$ again. The set inclusion diagrams of $\mathcal{M}_G$ and
$\mathcal{M}_{G*\{p,q\}}$ are given in
Figure~\ref{fig:ex_set_inclusion}.
\end{Example}

\paragraph{Graphs}
The pivots $G*X$ where $X$ is a minimal element of $\mathcal{M}_G
\backslash \{\emptyset\}$ w.r.t. inclusion are called
\emph{elementary}. It is noted in \cite{Geelen97} that an elementary
pivot $X$ corresponds to either a loop, $X = \{u\} \in E(G)$, or to
an edge, $X = \{u,v\} \in E(G)$, where both vertices $u$ and $v$ are
non-loops. Moreover, each $Y \in \mathcal{M}_G$ can be partitioned
$Y = X_1 \cup \cdots\cup X_n$ such that $G*Y = G*(X_1 \xor
\cdots\xor X_n) = (\cdots(G*X_1)\cdots * X_n)$ is a composition of
disjoint elementary pivots. Consequently, a direct definition of the
elementary pivots on graphs $G$ is sufficient to define the
(general) pivot operation.

The elementary pivot $G*\{u\}$ on a loop $\{u\}$ is called
\emph{local complementation}. It is the graph obtained from $G$ by
complementing the edges in the neighbourhood $N_G(u) = \{ v \in V
\mid \{u,v\} \in E(G), u \not= v \}$ of $u$ in $G$: for each $v,w
\in N_G(u)$, $\{v,w\}\in E(G)$ iff $\{v,w\} \not\in E(G*\{u\})$, and
$\{v\}\in E(G)$ iff $\{v\} \not\in E(G*\{u\})$ (the case $v=w$). The
other edges are left unchanged.

The elementary pivot $G*\pair uv$ on an edge $\pair uv$ between
distinct non-loop vertices $u$ and $v$ is called \emph{edge
complementation}. For a vertex $x$ consider its closed neighbourhood
$N'_G(x)= N_G(x)\cup \{x\}$. The edge $\pair uv$ partitions the
vertices of $G$ connected to $u$ or $v$ into three sets $V_1 =
N'_G(u) \setminus N'_G(v)$, $V_2 = N'_G(v) \setminus N'_G(u)$, $V_3
= N'_G(u) \cap N'_G(v)$. Note that $u,v \in V_3$.

\begin{figure}[t]
\centerline{\unitlength 1.0mm
\begin{picture}(55,42)(0,1)
\drawccurve(02,28)(25,21)(48,32)(25,39)
\drawccurve(00,10)(10,00)(20,10)(10,20)
\drawccurve(30,10)(40,00)(50,10)(40,20)
\gasset{AHnb=0,Nw=1.5,Nh=1.5,Nframe=n,Nfill=y}
\gasset{ExtNL=y,NLdist=1.5,NLangle=90}
\put(10,02){\makebox(0,0)[cc]{$V_1$}}
\put(40,02){\makebox(0,0)[cc]{$V_2$}}
\put(25,36){\makebox(0,0)[cc]{$V_3$}}
  \node(u)(09,28){$u$}
  \node(v)(20,30){$v$}
  \node(uu)(29,32){}
  \node(vv)(41,28){}
  \node(u1)(7,14){}
  \node(u2)(14,7){}
  \node(v1)(38,7){}
  \node(v2)(43,14){}
  \drawedge(u,v){}
  \drawedge(u,u1){}
  \drawedge(u,u2){}
  \drawedge(v,v1){}
  \drawedge(v,v2){}
  \drawedge(u1,v2){}
  \drawedge(u2,v1){}
  \drawedge[dash={1}0](v1,v2){}
  \drawedge[dash={1}0](u1,u2){}
  \drawedge[dash={1}0](uu,vv){}
  \drawedge(uu,u1){}
  \drawedge(vv,u2){}
  \drawedge(uu,v1){}
  \drawedge(vv,v2){}
\end{picture}
\begin{picture}(55,42)(0,1)
\drawccurve(02,28)(25,21)(48,32)(25,39)
\drawccurve(00,10)(10,00)(20,10)(10,20)
\drawccurve(30,10)(40,00)(50,10)(40,20)
\gasset{AHnb=0,Nw=1.5,Nh=1.5,Nframe=n,Nfill=y}
\gasset{ExtNL=y,NLdist=1.5,NLangle=90}
\put(10,02){\makebox(0,0)[cc]{$V_1$}}
\put(40,02){\makebox(0,0)[cc]{$V_2$}}
\put(25,36){\makebox(0,0)[cc]{$V_3$}}
  \node(u)(09,28){$u$}
  \node(v)(20,30){$v$}
  \node(uu)(29,32){}
  \node(vv)(41,28){}
  \node(u1)(7,14){}
  \node(u2)(14,7){}
  \node(v1)(38,7){}
  \node(v2)(43,14){}
  \drawedge(u,v){}
  \drawedge(v,u1){}
  \drawedge(v,u2){}
  \drawedge(u,v1){}
  \drawedge(u,v2){}
  \drawedge(u1,v1){}
  \drawedge(u2,v2){}
  \drawedge[dash={1}0](v1,v2){}
  \drawedge[dash={1}0](u1,u2){}
  \drawedge[dash={1}0](uu,vv){}
  \drawedge(uu,u2){}
  \drawedge(vv,u1){}
  \drawedge(uu,v2){}
  \drawedge(vv,v1){}
\end{picture}%
}
\caption{Pivoting $\pair uv$ in a graph.
 Connection $\pair xy$ is toggled iff $x\in V_i$ and $y\in V_j$ with
 $i\neq j$. Note that $u$ and $v$ are connected to all vertices in $V_3$,
 these edges are omitted in the diagram. The operation does not
 affect edges adjacent to vertices outside the sets
 $V_1,V_2,V_3$, nor does it change any of the loops.}%
\label{fig:pivot}
\end{figure}
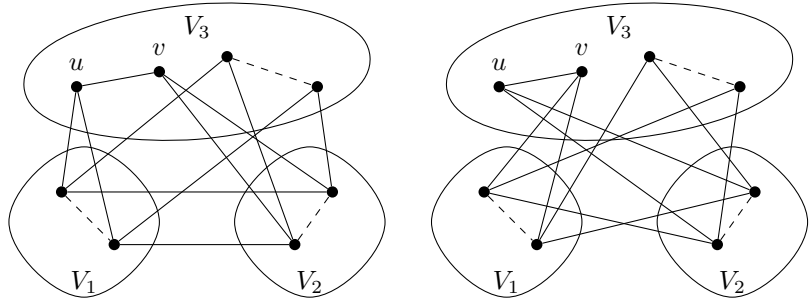

The graph  $G*\pair uv$ is constructed by ``toggling'' all edges
between different $V_i$ and $V_j$: for $\pair xy$ with $x\in V_i$
and $y\in V_j$ ($i\neq j$): $\pair xy \in E(G)$ iff $\pair xy \notin
E(G\sub{\{u,v\}})$, see Figure~\ref{fig:pivot}. The remaining edges
remain unchanged. Note that, as a result of this operation, the
neighbours of $u$ and $v$ are interchanged.

\begin{Example}
The whole orbit of $G$ of Example~\ref{ex:introd_pivot} under pivot
is given in Figure~\ref{fig:pivot_space}. It is obtained by
iteratively applying elementary pivots to $G$. Note that $G*\{p,q\}$
\emph{is} defined (top-right) but it is not an elementary pivot.
\end{Example}

\section{Dual Pivots}
In this section we introduce the dual pivot and show that it has
some interesting properties.

First note that the next result follows directly from
Equation~(\ref{pivot_def_reverse}).
\begin{Lemma} \label{lem:equaleigenspace}
Let $A$ be a $V \times V$-matrix (over some field) and let
$X \subseteq V$ with $A\sub{X}$ nonsingular.
Then the eigenspaces of $A$ and $A*X$ on value $1$ are equal, i.e.,
$E_1(A) = E_1(A*X)$.
\end{Lemma}
\begin{Proof}
We have $v \in E_1(A)$ iff $Av = v$ iff $(A*X)v = v$ iff $v \in E_1(A*X)$.
\end{Proof}

For a graph $G$, we denote $G+I$ to be the graph having adjacency
matrix $A(G)+I$ where $I$ is the identity matrix. Thus, $G+I$ is
obtained from $G$ by replacing each loop by a non-loop and vice
versa.

\begin{Definition}
Let $G$ be a graph and let $X \subseteq V$ with $\det ((G+I)\sub{X})
= 1$. The \emph{dual pivot} of $G$ on $X$, denoted by $G \bar* X$,
is $((G+I)*X)+I$.
\end{Definition}
Note that the condition $\det ((G+I)\sub{X}) = 1$ in the definition
of dual pivot ensures that the expression $((G+I)*X)+I$ is defined.
The dual pivot may be considered as the pivot operation
\emph{conjugated} by addition of the identity matrix $I$.
As $I+I$ is the null matrix (over $\two$), we have, similar as
for pivot, that dual pivot is an involution, and more
generally $(G\bar{*}X)\bar{*}Y$, when defined, is equal to
$G\bar{*}(X \xor Y)$.

By Lemma~\ref{lem:equaleigenspace}, we have the following result.
\begin{Lemma} \label{lem:dualp_equal_kernel}
Let $G$ be a graph and let $X \subseteq V$ such that $G\bar{*}X$ is
defined. Then $\kernel(G\bar{*}X) = \kernel(G)$.
\end{Lemma}
\begin{Proof}
Note that $Ax = 0$ iff $(A+I)x = x$. Hence, $\kernel(G) = E_1(G+I)$.
Since $I+I$ is the null matrix (over $\two$), we have also
$\kernel(G+I) = E_1(G)$.

Therefore, we have $\kernel(G\bar{*}X) = \kernel(((G+I)*X)+I) =
E_1((G+I)*X) = E_1(G+I)$, where we used
Lemma~\ref{lem:equaleigenspace} is the last equality. Finally,
$E_1(G+I) = \kernel(G)$ and therefore we obtain $\kernel(G\bar{*}X)
= \kernel(G)$.
\end{Proof}
In particular, for the case $X = V$, we have that
$\kernel((G+I)^{-1}+I) = \kernel(G)$ (the inverse is computed over
$\two$) if the left-hand side is defined.

\begin{Remark}
By Lemma~\ref{lem:dualp_equal_kernel} and
Corollary~\ref{cor:ker_iff_max_indep} we have that the (column)
matroids associated with $G$ and $G\bar{*}X$ are equal. Note that
here the matroids are obtained from the column vectors of the
adjacency matrices of $G$ and $G\bar{*}X$; this is not to be
confused with graphic matroids which are obtained from the column
vectors of the incidence matrices of graphs.
\end{Remark}

We call dual pivot $G \bar* X$ \emph{elementary} if $*X$ is an
elementary pivot for $G+I$. Equivalently, they are the dual pivots
on $X$ for which there is no non-empty $Y \subset X$ where
$G\bar{*}Y$ is applicable. An elementary dual pivot $\bar*\{u\}$ is
defined on a non-loop vertex $u$, and an elementary dual pivot
$\bar*\{u,v\}$ is defined on an edge $\{u,v\}$ where both $u$ and
$v$ have loops. This is the only difference between pivot and its
dual: both the elementary dual pivot $\bar*\{u\}$ and the elementary
pivot $*\{u\}$ have the same effect on the graph --- both ``take the
complement'' of the neighbourhood of $u$. Similarly, the effect of
the elementary dual pivot $\bar*\{u,v\}$ and the elementary pivot
$*\{u,v\}$ is the same, only the condition when they can be applied
differs.

Note that the eigenspaces $E_0(G) = \kernel(G)$ and $E_1(G)$ have a
natural interpretation in graph terminology. For $X \subseteq V(G)$,
$X \in E_0(G)$ iff every vertex in $V(G)$ is connected to an even
number of vertices in $X$ (loops do count). Also, $X \in E_1(G)$ iff
every vertex in $V(G) \backslash X$ is connected to an even number
of vertices in $X$ and every vertex in $X$ is connected to an odd
number of vertices in $X$ (again loops do count).

\begin{figure}[t]
\unitlength 0.7mm%
\begin{center}
\begin{picture}(100,80)
\put(0,50){
\begin{picture}(30,30)
\gasset{AHnb=0,Nw=1.5,Nh=1.5,Nframe=n,Nfill=y}
\gasset{AHnb=0,Nw=8,Nh=8,Nframe=y,Nfill=n}
  \node(5)(05,05){$r$}
  \node(2)(05,25){$p$}
  \node(3)(25,25){$q$}
  \node(6)(25,05){$s$}
  \drawedge(2,3){}
  \drawedge(2,5){}
  \drawedge(3,6){}
  \drawedge(5,6){}
  \drawedge(2,6){}
  \drawloop[loopangle=135,loopdiam=7](2){}
  \drawloop[loopangle=-135,loopdiam=7](5){}
  \drawloop[loopangle=-45,loopdiam=7](6){}
\end{picture}}
\node[Nframe=n](1)(35,65){} \node[Nframe=n](2)(70,65){}
\drawedge[AHnb=1,ATnb=1](1,2){$*\{p\}$}
\node[Nframe=n](3)(35,15){} \node[Nframe=n](4)(70,15){}
\drawedge[AHnb=1,ATnb=1](3,4){$\bar{*}\{p\}$}
\node[Nframe=n](5)(17,27){} \node[Nframe=n](6)(17,53){}
\drawedge[AHnb=1,ATnb=1](5,6){$+I$}
\node[Nframe=n](5)(87,27){} \node[Nframe=n](6)(87,53){}
\drawedge[AHnb=1,ATnb=1](5,6){$+I$}
\put(70,50){
\begin{picture}(30,30)
\gasset{AHnb=0,Nw=1.5,Nh=1.5,Nframe=n,Nfill=y}
\gasset{AHnb=0,Nw=8,Nh=8,Nframe=y,Nfill=n}
  \node(5)(05,05){$r$}
  \node(2)(05,25){$p$}
  \node(3)(25,25){$q$}
  \node(6)(25,05){$s$}
  \drawedge(2,3){}
  \drawedge(2,5){}
  \drawedge(3,5){}
  \drawedge(2,6){}
  \drawloop[loopangle=135,loopdiam=7](2){}
  \drawloop[loopangle=45,loopdiam=7](3){}
\end{picture}}
\put(0,0){
\begin{picture}(30,30)
\gasset{AHnb=0,Nw=1.5,Nh=1.5,Nframe=n,Nfill=y}
\gasset{AHnb=0,Nw=8,Nh=8,Nframe=y,Nfill=n}
  \node(5)(05,05){$r$}
  \node(2)(05,25){$p$}
  \node(3)(25,25){$q$}
  \node(6)(25,05){$s$}
  \drawedge(2,3){}
  \drawedge(2,5){}
  \drawedge(3,6){}
  \drawedge(5,6){}
  \drawedge(2,6){}
  \drawloop[loopangle=45,loopdiam=7](3){}
\end{picture}}
\put(70,0){
\begin{picture}(30,30)
\gasset{AHnb=0,Nw=1.5,Nh=1.5,Nframe=n,Nfill=y}
\gasset{AHnb=0,Nw=8,Nh=8,Nframe=y,Nfill=n}
  \node(5)(05,05){$r$}
  \node(2)(05,25){$p$}
  \node(3)(25,25){$q$}
  \node(6)(25,05){$s$}
  \drawedge(2,3){}
  \drawedge(2,5){}
  \drawedge(3,5){}
  \drawedge(2,6){}
  \drawloop[loopangle=-135,loopdiam=7](5){}
  \drawloop[loopangle=-45,loopdiam=7](6){}
\end{picture}}
\end{picture}
\end{center}
\caption{Dual pivot of graph $G$ from Example~\ref{ex:introd_pivot}
($G$ is shown in the lower-left corner).} \label{fig:c4}
\end{figure}
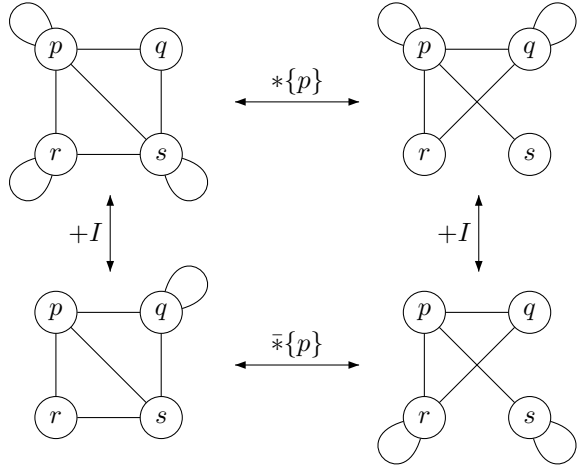

\begin{Example}\label{ex:c4}
Let $G'$ be the graph depicted on the upper-left corner of
Figure~\ref{fig:c4}. We have $A(G') = \kbordermatrix{
  & p & q & r & s \\
p & 1 & 1 & 1 & 1 \\
q & 1 & 0 & 0 & 1 \\
r & 1 & 0 & 1 & 1 \\
s & 1 & 1 & 1 & 1 }$. Note that $E_1(G') = \{ \varnothing, \{p,r,s\}
\}$ is of dimension $1$. We can apply an elementary pivot over $p$
on $G'$. The resulting graph $G'*\{p\}$ is depicted on the
upper-right corner of Figure~\ref{fig:c4}, and we have $A(G'*\{p\})
=
\kbordermatrix{
  & p & q & r & s \\
p & 1 & 1 & 1 & 1 \\
q & 1 & 1 & 1 & 0 \\
r & 1 & 1 & 0 & 0 \\
s & 1 & 0 & 0 & 0 }
.$
Note that the elements of $E_1(G')$ are precisely the eigenvectors
(or eigensets) on $1$ for $A(G'*\{p\})$, cf.
Lemma~\ref{lem:equaleigenspace}. The graphs $G'+I$ (which is $G$ in
Example~\ref{ex:introd_pivot}) and $G'*\{p\}+I$ are depicted in the
lower-left and lower-right corner of Figure~\ref{fig:c4},
respectively. By definition of the dual pivot we have $G
\bar{*}\{p\} = (G'+I)\bar{*}\{p\} = G'*\{p\}+I$.
\end{Example}

It is a basic fact from linear algebra that elementary row
operations retain the kernel of matrices.
Lemma~\ref{lem:dualp_equal_kernel} suggests that the dual pivot may
possibly be simulated by elementary row operations. We now show that
this is indeed the case. Over $\two$ the elementary row operations
are 1) row switching and 2) adding one row to another (row
multiplication over $\two$ does not change the matrix). The
elementary row operations corresponding to the dual pivot operation
are easily deduced by restricting to elementary dual pivots. The
dual pivot on a non-loop vertex $u$ corresponds, in the adjacency
matrix, to adding the row corresponding to $u$ to each row
corresponding to a vertex in the neighbourhood of $u$. Moreover, the
dual pivot on edge $\{u,v\}$ (where both $u$ and $v$ have loops)
corresponds to 1) adding the row corresponding to $u$ to each row
corresponding to a vertex in the neighbourhood of $v$ except $u$, 2)
adding the row corresponding to $v$ to each row corresponding to a
vertex in the neighbourhood of $u$ except $v$, 3) switching the rows
of $u$ and $v$. Note that this procedure allows for another,
equivalent, definition of the regular pivot: add $I$, apply the
corresponding elementary row operations, and finally add $I$ again.

Note that the dual pivot has the property that it transforms a
symmetric matrix to another symmetric matrix with equal kernel.
Applying elementary row operations however will in general not
obtain symmetric matrices.

\section{Maximal Pivots}
In Section~\ref{sec:def_pivots} we recalled that the minimal
elements of $\mathcal{M}_G$, corresponding to elementary pivots,
form the building blocks of (general) pivots. In this section we
show that the set of maximal elements of $\mathcal{M}_G$,
corresponding to ``maximal pivots'', is invariant under dual pivot.

For $\mathcal{M}_G = (V,D_G)$, we define $\mathcal{F}_G =
\max(D_G)$. Thus, for $X \subseteq V(G)$, $X \in \mathcal{F}_G$ iff
$\det G\sub{X} = 1$ while $\det G\sub{Y} = 0$ for every $Y \supset
X$.

\begin{figure}[bt]
\unitlength 0.7mm
\begin{center}
\begin{picture}(50,80)(-5,-20)
\gasset{Nw=8,Nh=8,Nframe=y,Nfill=n}
  \node(nul)(15,60){$\varnothing$}
  \node(q)(15,40){$q$}
  \node(pr)(00,20){$pr$}
  \node(pq)(10,20){$pq$}
  \node(qs)(20,20){$qs$}
  \node(ps)(30,20){$ps$}
  \node(rs)(40,20){$rs$}
  \node(pqr)(05,00){$pqr$}
  \node(pqs)(20,00){$pqs$}
  \node(qrs)(35,00){$qrs$}
  \drawedge(nul,q){}
  \drawedge(nul,pr){}
  \drawedge(nul,ps){}
  \drawedge(nul,rs){}
  \drawedge(pr,pqr){}
  \drawedge(ps,pqs){}
  \drawedge(rs,qrs){}
  \drawedge(q,pq){}
  \drawedge(q,qs){}
  \drawedge(pq,pqr){}
  \drawedge(pq,pqs){}
  \drawedge(qs,pqs){}
  \drawedge(qs,qrs){}
\end{picture}
\begin{picture}(50,80)(-5,-20)
\gasset{Nw=8,Nh=8,Nframe=y,Nfill=n}
 \node(nul)(25,60){$\varnothing$}
 \node(r)(15,40){$r$}
 \node(p)(25,40){$p$}
 \node(s)(35,40){$s$}
 \node(qr)(00,20){$qr$}
 \node(pq)(10,20){$pq$}
 \node(ps)(30,20){$ps$}
 \node(qs)(40,20){$qs$}
 \node(prs)(22,00){$prs$}
 \node(pqrs)(20,-20){$pqrs$}
 \drawedge(nul,p){}
 \drawedge(nul,r){}
 \drawedge(nul,s){}
 \drawedge(p,pq){}
 \drawedge(p,ps){}
 \drawedge(r,qr){}
 \drawedge(s,ps){}
 \drawedge(s,qs){}
 \drawedge(r,prs){}
 \drawedge(ps,prs){}
 \drawedge(qs,pqrs){}
 \drawedge(pq,pqrs){}
 \drawedge(qr,pqrs){}
 \drawedge(prs,pqrs){}
\end{picture}
\begin{picture}(50,80)(-5,-20)
\gasset{Nw=8,Nh=8,Nframe=y,Nfill=n}
  \node(nul)(15,60){$\varnothing$}
  \node(r)(05,40){$r$}
  \node(s)(20,40){$s$}
  \node(pr)(00,20){$pr$}
  \node(qr)(10,20){$qr$}
  \node(rs)(20,20){$rs$}
  \node(ps)(30,20){$ps$}
  \node(pq)(40,20){$pq$}
  \node(pqr)(05,00){$pqr$}
  \node(pqs)(35,00){$pqs$}
  \node(qrs)(20,00){$qrs$}
  \drawedge(nul,r){}
  \drawedge(nul,s){}
  \drawedge(nul,pq){}
  \drawedge(pr,pqr){}
  \drawedge(ps,pqs){}
  \drawedge(rs,qrs){}
  \drawedge(qr,pqr){}
  \drawedge(qr,qrs){}
  \drawedge(pq,pqr){}
  \drawedge(pq,pqs){}
  \drawedge(r,pr){}
  \drawedge(r,qr){}
  \drawedge(r,rs){}
  \drawedge(s,rs){}
  \drawedge(s,ps){}
\end{picture}
\end{center}
\caption{Set inclusion diagram of $\mathcal{M}_G$,
$\mathcal{M}_{G{*}\{p,q\}}$, and $\mathcal{M}_{G\bar{*}\{p\}}$ for
$G$, $G{*}\{p,q\}$, and $G\bar{*}\{p\}$ as given in
Examples~\ref{ex:introd_pivot} and \ref{ex:c4}.}
\label{fig:ex_set_inclusion}
\end{figure}
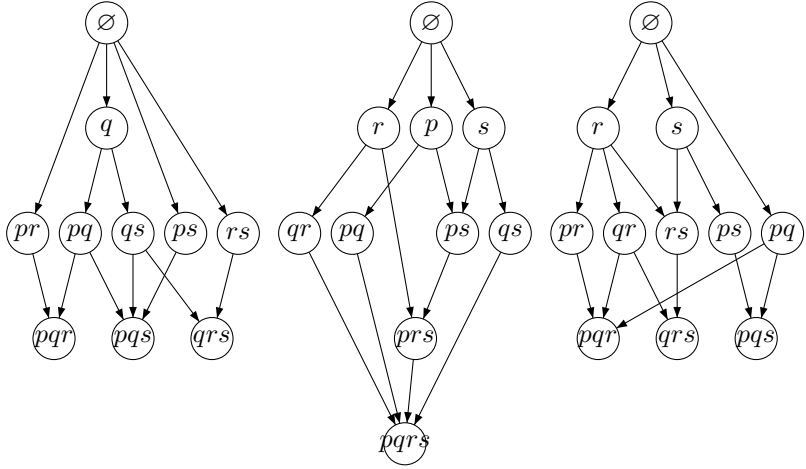

\begin{Example}
We continue Example~\ref{ex:introd_pivot}. Let $G$ be the graph on
the lower-left corner of Figure~\ref{fig:c4}. Then from the set
inclusion diagram of $\mathcal{M}_G$ in
Figure~\ref{fig:ex_set_inclusion} we see that $\mathcal{F}_{G} =
\{\{p,q,s\},\{p,q,r\},\{q,r,s\}\}$. Also we see from the figure that
$\mathcal{F}_{G*\{p,q\}} = \{V\}$.
\end{Example}

Next we recall the Strong Principal Minor Theorem for
\mbox{(quasi-)} symmetric matrices from
\cite{Kodiyalam_Lam_Swan_2008}
--- it is stated here for graphs (i.e., symmetric matrices over $\two$).
\footnote{Clearly, for a matrix $A$, $\det A\sub{X} \not= 0$ implies
that $X$ is independent for $A$. The reverse implication is not
valid in general.}

\begin{Proposition} \label{prop:spmt}
Let $G$ be a graph such that $A(G)$ has rank $r$, and let $X
\subseteq V(G)$ with $|X| = r$. Then $X$ is independent for $A(G)$
iff $\det G\sub{X} = 1$.
\end{Proposition}

Note that the independent sets $X$ of cardinality equal to the rank
are precisely the bases of a matrix $A$.

The following result is easy to see now from
Proposition~\ref{prop:spmt}.
\begin{Lemma} \label{lem:all_in_F_size_r}
Let $G$ be a graph such that $A(G)$ has rank $r$. Each element of
$\mathcal{F}_G$ is of cardinality $r$.
\end{Lemma}
\begin{Proof}
If there is an $X \in \mathcal{F}_G$ of cardinality $q > r$, then
the columns of $A(G\sub{X})$ are linearly independent, and thus so
are the columns of $A(G)$ corresponding to $X$. This contradicts the
rank of $A(G)$.

Finally, assume that there is an $X \in \mathcal{F}_G$ of
cardinality $q < r$. Since the columns of $A(G\sub{X})$ are linearly
independent, so are the columns of $A(G)$ corresponding to $X$.
Since $A(G)$ has rank $r$, $X$ can be extended to a set $X'$ with
cardinality $r$. Hence by Proposition~\ref{prop:spmt} $\det
G\sub{X'} = 1$ with $X' \supset X$ --- a contradiction of $X \in
\mathcal{F}_G$.
\end{Proof}

\begin{Example}
We continue Example~\ref{ex:introd_pivot}. Let again $G$ be the
graph on the lower-left corner of Figure~\ref{fig:c4}. Then the
elements $\mathcal{F}_{G} = \{\{p,q,s\},\{p,q,r\},\{q,r,s\}\}$ are
all of cardinality $3$ --- the rank of $A(G)$. Moreover,
$\mathcal{F}_{G*\{p,q\}} = \{V\}$ and $|V| = 4$ is equal to the rank
of $G*\{p,q\}$.
\end{Example}

Combining Proposition~\ref{prop:spmt} and
Lemma~\ref{lem:all_in_F_size_r}, we have the following result.
\begin{Corollary} \label{cor:char_extremal_basis}
Let $G$ be a graph, and let $X \subseteq V(G)$. Then $X$ is a basis
for $A(G)$ iff $X \in \mathcal{F}_G$.
\end{Corollary}
Equivalently, with $\mathcal{P}_G = (I,D)$ from
Section~\ref{sec:notation}, Corollary~\ref{cor:char_extremal_basis}
states that $\max(I) = \mathcal{F}_G$.

By Corollaries~\ref{cor:ker_iff_max_indep} and
\ref{cor:char_extremal_basis} we have now the following.
\begin{Lemma} \label{lem:equal_kernel_maxminors}
Let $G$ and $G'$ be graphs. Then $\mathcal{F}_{G} =
\mathcal{F}_{G'}$ iff $\kernel(G) = \kernel(G')$.
\end{Lemma}

Recall that Lemma~\ref{lem:dualp_equal_kernel} shows that the dual
pivot retains the kernel. We may now conclude from
Lemma~\ref{lem:equal_kernel_maxminors} that also $\mathcal{F}_G$ is
retained under dual pivot. It is the main result of this paper, and,
as we will see in Section~\ref{sec:appl_ga}, has an important
application.
\begin{Theorem} \label{thm:dual}
Let $G$ be a graph, and let $X \subseteq V$. Then $\mathcal{F}_{G} =
\mathcal{F}_{G \bar* X}$ if the right-hand side is defined.
\end{Theorem}

In particular, the case $X=V$, we have $\mathcal{F}_{G+I} =
\mathcal{F}_{G^{-1}+I}$ if $G$ is invertible (over $\two$).

Let $\mathcal{O}_G = \{ G\bar{*}X \mid X \subseteq V, \det
(G+I)\sub{X} = 1 \}$ be the orbit of $G$ under dual pivot, and note
that $G \in \mathcal{O}_G$. By Theorem~\ref{thm:dual}, if $G_1, G_2
\in \mathcal{O}_G$, then $\mathcal{F}_{G_1} = \mathcal{F}_{G_2}$.
Note that the reverse implication does not hold: e.g. $\mathcal{O}_I
= \{I\}$ and $\mathcal{F}_I = \{V\}$, while clearly there are many
other graphs $G$ with $\det G = 1$ (which means $\mathcal{F}_G =
\{V\}$).

\begin{Example} \label{ex:equal_max_F_dual_pivot}
We continue Example~\ref{ex:c4}. Let again $G$ be the graph on the
lower-left corner of Figure~\ref{fig:c4}. Then $G\bar{*}\{p\}$ is
depicted on the lower-right corner of Figure~\ref{fig:c4}. We have
$\mathcal{F}_{G \bar* \{p\}} = \{\{p,q,s\},\{p,q,r\},\{q,r,s\}\}$,
see Figure~\ref{fig:ex_set_inclusion}, so indeed $\mathcal{F}_{G} =
\mathcal{F}_{G \bar*\{p\}}$.
\end{Example}

For symmetric $V \times V$-matrices $A$ over $\two$,
Theorem~\ref{thm:dual} states that if $A$ can be partially inverted
w.r.t. $Y \subseteq V$, where $Y$ is maximal w.r.t. set inclusion,
then this holds for every matrix obtained from $A$ by dual pivot.

\section{Maximal Contractions}
For a graph $G$, we define the \emph{contraction} of $G$ on $X
\subseteq V$ with $\det G\sub{X} = 1$, denoted by $G *\backslash X$,
to be the graph $(G*X) \backslash X$ --- the pivot on $X$ followed
by the removal of the vertices of $X$. Equivalently, contraction is
the Schur complement applied to graphs. A contraction of $G$ on $X$
is \emph{maximal} if there is no $Y \supset X$ such that $\det
G\sub{Y} = 1$, hence if $X \in \mathcal{F}_{G}$. The graph obtained
by a maximal contraction on $X$ is a discrete graph $G'$ (without
loops). Indeed, if $G'$ were to have a loop $e = \{u\}$ or an edge
$e = \{u,v\}$ between two non-loop vertices, then, since $\det
G\sub{X \xor e} = \det ((G*X)\sub{e}) = \det ((G*\backslash
X)\sub{e}) = 1$, $X \xor e \supset X$ would be a contradiction of
the maximality of $X$.
Moreover, by Lemma~\ref{lem:all_in_F_size_r}, the number of vertices
of $G'$ is equal to the nullity (dimension of the kernel, which
equals the dimension of the matrix minus its rank) of $G$.

\begin{Remark}
In fact, it is known that any Schur complement in a matrix $A$ has
the same nullity as $A$ itself --- it is a consequence of the
Guttman rank additivity formula, see, e.g.,
\cite[Section~6.0.1]{SchurBook2005}. Therefore the the nullity is
invariant under contraction in general (not only maximal
contraction).
\end{Remark}


By Theorem~\ref{thm:dual} we have the following.
\begin{Corollary} \label{cor:max_contractions}
The set of discrete graphs obtainable through contractions is equal
for $G$ and $G\bar{*}X$ for all $X \subseteq V$ with $\det
(G+I)\sub{X} = 1$.
\end{Corollary}
In this sense, all the elements of the orbit $\mathcal{O}_G$ have
equal ``behaviour'' w.r.t. maximal contractions.

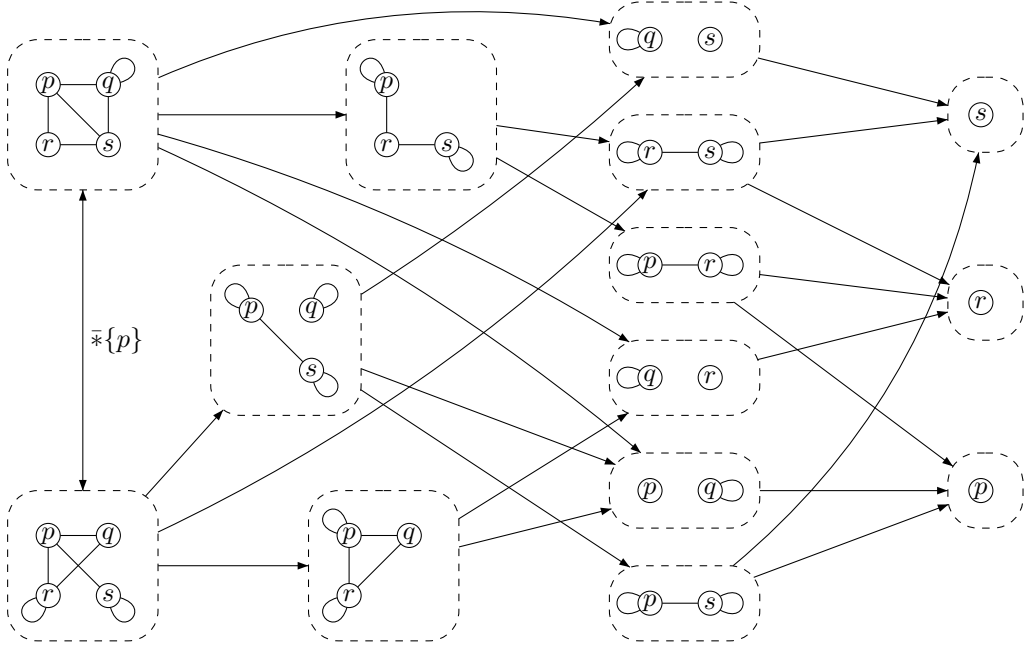
\begin{figure}[t]
\unitlength 1mm%
\begin{center}
\begin{picture}(120,80)(00,-05)
\node[Nw=20,Nh=20,Nframe=y,dash={1}0](I)(00,60){}
\node[Nw=20,Nh=20,Nframe=n](I)(00,60){%
  \unitlength0.4mm
  \begin{picture}(30,30)
  \gasset{AHnb=0,Nw=1.5,Nh=1.5,Nframe=n,Nfill=y}
  \gasset{AHnb=0,Nw=8,Nh=8,Nframe=y,Nfill=n}
    \node(r)(05,05){$r$}
    \node(p)(05,25){$p$}
    \node(q)(25,25){$q$}
    \node(s)(25,05){$s$}
    \drawedge(p,q){}
    \drawedge(p,r){}
    \drawedge(p,s){}
    \drawedge(q,s){}
    \drawedge(r,s){}
    \drawloop[loopangle=45,loopdiam=7](q){}
  \end{picture}
}
\node[Nw=20,Nh=20,Nframe=y,dash={1}0](prs)(45,60){}
\node[Nw=20,Nh=20,Nframe=n](prs)(45,60){%
  \unitlength0.4mm
  \begin{picture}(30,30)
  \gasset{AHnb=0,Nw=1.5,Nh=1.5,Nframe=n,Nfill=y}
  \gasset{AHnb=0,Nw=8,Nh=8,Nframe=y,Nfill=n}
  \node(r)(05,05){$r$}
  \node(p)(05,25){$p$}
  \node(s)(25,05){$s$}
  \drawedge(p,r){}
  \drawedge(r,s){}
  \drawloop[loopangle=135,loopdiam=7](p){}
  \drawloop[loopangle=-45,loopdiam=7](s){}
\end{picture}
}
\node[Nw=20,Nh=20,Nframe=y,dash={1}0](pqs)(27,30){}
\node[Nw=20,Nh=20,Nframe=n](pqs)(27,30){%
  \unitlength0.4mm
  \begin{picture}(30,30)
  \gasset{AHnb=0,Nw=1.5,Nh=1.5,Nframe=n,Nfill=y}
  \gasset{AHnb=0,Nw=8,Nh=8,Nframe=y,Nfill=n}
  \node(p)(05,25){$p$}
  \node(q)(25,25){$q$}
  \node(s)(25,05){$s$}
  \drawedge(p,s){}
  \drawloop[loopangle=135,loopdiam=7](p){}
  \drawloop[loopangle=45,loopdiam=7](q){}
  \drawloop[loopangle=-45,loopdiam=7](s){}
\end{picture}
}
\node[Nw=20,Nh=20,Nframe=y,dash={1}0](pqr)(40,00){}
\node[Nw=20,Nh=20,Nframe=n](pqr)(40,00){%
  \unitlength0.4mm
  \begin{picture}(30,30)
  \gasset{AHnb=0,Nw=1.5,Nh=1.5,Nframe=n,Nfill=y}
  \gasset{AHnb=0,Nw=8,Nh=8,Nframe=y,Nfill=n}
  \node(r)(05,05){$r$}
  \node(p)(05,25){$p$}
  \node(q)(25,25){$q$}
  \drawedge(p,q){}
  \drawedge(p,r){}
  \drawedge(q,r){}
  \drawloop[loopangle=135,loopdiam=7](p){}
  \drawloop[loopangle=225,loopdiam=7](r){}
\end{picture}
}
\node[Nw=20,Nh=10,Nframe=y,dash={1}0](qs)(80,70){}
\node[Nw=20,Nh=10,Nframe=n](qs)(80,70){%
  \unitlength0.4mm
  \begin{picture}(30,10)
  \gasset{AHnb=0,Nw=1.5,Nh=1.5,Nframe=n,Nfill=y}
  \gasset{AHnb=0,Nw=8,Nh=8,Nframe=y,Nfill=n}
  \node(1)(05,05){$q$}
  \node(2)(25,05){$s$}
  \drawloop[loopangle=180,loopdiam=7](1){}
\end{picture}
}
\node[Nw=20,Nh=10,Nframe=y,dash={1}0](rs)(80,55){}
\node[Nw=20,Nh=10,Nframe=n](rs)(80,55){%
  \unitlength0.4mm
  \begin{picture}(30,10)
  \gasset{AHnb=0,Nw=1.5,Nh=1.5,Nframe=n,Nfill=y}
  \gasset{AHnb=0,Nw=8,Nh=8,Nframe=y,Nfill=n}
  \node(1)(05,05){$r$}
  \node(2)(25,05){$s$}
  \drawedge(1,2){}
  \drawloop[loopangle=180,loopdiam=7](1){}
  \drawloop[loopangle=0,loopdiam=7](2){}
\end{picture}
}
\node[Nw=20,Nh=10,Nframe=y,dash={1}0](pr)(80,40){}
\node[Nw=20,Nh=10,Nframe=n](pr)(80,40){%
  \unitlength0.4mm
  \begin{picture}(30,10)
  \gasset{AHnb=0,Nw=1.5,Nh=1.5,Nframe=n,Nfill=y}
  \gasset{AHnb=0,Nw=8,Nh=8,Nframe=y,Nfill=n}
  \node(1)(05,05){$p$}
  \node(2)(25,05){$r$}
  \drawedge(1,2){}
  \drawloop[loopangle=180,loopdiam=7](1){}
  \drawloop[loopangle=0,loopdiam=7](2){}
\end{picture}
}
\node[Nw=20,Nh=10,Nframe=y,dash={1}0](qr)(80,25){}
\node[Nw=20,Nh=10,Nframe=n](qr)(80,25){%
  \unitlength0.4mm
  \begin{picture}(30,10)
  \gasset{AHnb=0,Nw=1.5,Nh=1.5,Nframe=n,Nfill=y}
  \gasset{AHnb=0,Nw=8,Nh=8,Nframe=y,Nfill=n}
  \node(1)(05,05){$q$}
  \node(2)(25,05){$r$}
  \drawloop[loopangle=180,loopdiam=7](1){}
\end{picture}
}
\node[Nw=20,Nh=10,Nframe=y,dash={1}0](pq)(80,10){}
\node[Nw=20,Nh=10,Nframe=n](pq)(80,10){%
  \unitlength0.4mm
  \begin{picture}(30,10)
  \gasset{AHnb=0,Nw=1.5,Nh=1.5,Nframe=n,Nfill=y}
  \gasset{AHnb=0,Nw=8,Nh=8,Nframe=y,Nfill=n}
  \node(1)(05,05){$p$}
  \node(2)(25,05){$q$}
  \drawloop[loopangle=0,loopdiam=7](2){}
\end{picture}
}
\node[Nw=20,Nh=10,Nframe=y,dash={1}0](ps)(80,-05){}
\node[Nw=20,Nh=10,Nframe=n](ps)(80,-05){%
  \unitlength0.4mm
  \begin{picture}(30,10)
  \gasset{AHnb=0,Nw=1.5,Nh=1.5,Nframe=n,Nfill=y}
  \gasset{AHnb=0,Nw=8,Nh=8,Nframe=y,Nfill=n}
  \node(1)(05,05){$p$}
  \node(2)(25,05){$s$}
  \drawedge(1,2){}
  \drawloop[loopangle=180,loopdiam=7](1){}
  \drawloop[loopangle=0,loopdiam=7](2){}
\end{picture}
}
\node[Nw=10,Nh=10,Nframe=y,dash={1}0](s1)(120,60){}
\node[Nw=10,Nh=10,Nframe=n](s1)(120,60){%
  \unitlength0.4mm
  \begin{picture}(10,10)
  \gasset{AHnb=0,Nw=1.5,Nh=1.5,Nframe=n,Nfill=y}
  \gasset{AHnb=0,Nw=8,Nh=8,Nframe=y,Nfill=n}
  \node(1)(05,05){$s$}
\end{picture}
}
\node[Nw=10,Nh=10,Nframe=y,dash={1}0](r1)(120,35){}
\node[Nw=10,Nh=10,Nframe=n](r1)(120,35){%
  \unitlength0.4mm
  \begin{picture}(10,10)
  \gasset{AHnb=0,Nw=1.5,Nh=1.5,Nframe=n,Nfill=y}
  \gasset{AHnb=0,Nw=8,Nh=8,Nframe=y,Nfill=n}
  \node(1)(05,05){$r$}
\end{picture}
}
\node[Nw=10,Nh=10,Nframe=y,dash={1}0](p1)(120,10){}
\node[Nw=10,Nh=10,Nframe=n](p1)(120,10){%
  \unitlength0.4mm
  \begin{picture}(10,10)
  \gasset{AHnb=0,Nw=1.5,Nh=1.5,Nframe=n,Nfill=y}
  \gasset{AHnb=0,Nw=8,Nh=8,Nframe=y,Nfill=n}
  \node(1)(05,05){$p$}
\end{picture}
}
\node[Nw=20,Nh=20,Nframe=y,dash={1}0](V)(00,00){}
\node[Nw=20,Nh=20,Nframe=n](V)(00,00){%
  \unitlength0.4mm
  \begin{picture}(30,30)
  \gasset{AHnb=0,Nw=1.5,Nh=1.5,Nframe=n,Nfill=y}
  \gasset{AHnb=0,Nw=8,Nh=8,Nframe=y,Nfill=n}
  \node(r)(05,05){$r$}
  \node(p)(05,25){$p$}
  \node(q)(25,25){$q$}
  \node(s)(25,05){$s$}
  \drawedge(p,q){}
  \drawedge(p,r){}
  \drawedge(p,s){}
  \drawedge(q,r){}
  \drawloop[loopangle=225,loopdiam=7](r){}
  \drawloop[loopangle=-45,loopdiam=7](s){}
\end{picture}
}
  \drawedge[AHnb=1](I,prs){}
  \drawedge[AHnb=1](V,pqs){}
  \drawedge[AHnb=1](V,pqr){}
  \drawedge[AHnb=1,curvedepth=8](I,qs){}
  \drawedge[AHnb=1](prs,rs){}
  \drawedge[AHnb=1](prs,pr){}
  \drawedge[AHnb=1,curvedepth=4](I,qr){}
  \drawedge[AHnb=1,curvedepth=4](I,pq){}
  \drawedge[AHnb=1,curvedepth=-5](V,rs){}
  \drawedge[AHnb=1,curvedepth=-2](pqs,qs){}
  \drawedge[AHnb=1](pqs,pq){}
  \drawedge[AHnb=1](pqs,ps){}
  \drawedge[AHnb=1](pqr,qr){}
  \drawedge[AHnb=1](pqr,pq){}
  \drawedge[AHnb=1](qs,s1){}
  \drawedge[AHnb=1](rs,s1){}
  \drawedge[AHnb=1](rs,r1){}
  \drawedge[AHnb=1](pr,r1){}
  \drawedge[AHnb=1](pr,p1){}
  \drawedge[AHnb=1](qr,r1){}
  \drawedge[AHnb=1](pq,p1){}
  \drawedge[AHnb=1](ps,p1){}
  \drawedge[AHnb=1,curvedepth=-8](ps,s1){}
%
%
%
%
%
  \drawedge[AHnb=1,ATnb=1](I,V){$\bar{*}\{p\}$}
\end{picture}
\end{center}
\caption{Elementary contractions starting from $G$ and
$G\bar{*}\{p\}$.} \label{fig:ex_max_contractions}
\end{figure}

\begin{Example} \label{ex:max_contractions}
We continue the example. Recall that, from
Example~\ref{ex:equal_max_F_dual_pivot}, $\mathcal{F}_{G} =
\mathcal{F}_{G\bar{*}\{p\}} = \{\{p,q,s\},\{p,q,r\},\{q,r,s\}\}$.
The elementary contractions starting from $G$ and $G\bar{*}\{p\}$
are given in Figure~\ref{fig:ex_max_contractions}. Notice that the
maximal contractions of $G$ and $G\bar{*}\{p\}$ obtain the same set
of (discrete) graphs.
\end{Example}

It is important to realize that while the \emph{maximal}
contractions (corresponding to $\mathcal{F}_{G}$) are the same for
graphs $G$ and $G\bar{*}X$, the whole set of contractions
(corresponding to $\mathcal{M}_G$) may be spectacularly different.
Indeed, e.g., in Example~\ref{ex:max_contractions}, the elementary
pivots for $G$ are $*\{q\}$, $*\{p,s\}$, $*\{p,r\}$, and $*\{r,s\}$,
while the elementary pivots for $G\bar{*}\{p\}$ are $*\{r\}$,
$*\{s\}$, and $*\{p,q\}$ (see Figure~\ref{fig:ex_set_inclusion}).

\section{Application: Gene Assembly} \label{sec:appl_ga}

Gene assembly is a highly involved and parallel process occurring in
one-cellular organisms called ciliates. During gene assembly a
nucleus, called micronucleus (MIC), is transformed into another
nucleus called macronucleus (MAC). Segments of the genes in the MAC
occur in scrambled order in the MIC \cite{GeneAssemblyBook}. During
gene assembly, recombination takes place to ``sort'' these gene
segments in the MIC in the right orientation and order to obtain the
MAC gene. The transformation of single genes from their MIC form to
their MAC form is formally modelled, see
\cite{Equiv_String_Graph_1,Equiv_String_Graph_2,GeneAssemblyBook},
as both a string based model and a (almost equivalent) graph based
model. It is observed in \cite{BHH/PivotsDetPM/09} that two of the
three operations in the graph based model are exactly the two
elementary principal pivot transform (PPT, or simply \emph{pivot})
operations on the corresponding adjacency matrices considered over
$\two$. The third operation simply removes isolated vertices.

Maximal contractions are especially important within the theory of
gene assembly in ciliates --- such a maximal sequence determines a
complete transformation of the gene to its MAC form. We first recall
the string rewriting system, and then recall the generalization to
the graph rewriting system.

\newcommand{\pset}[1]{\|#1\|}

Let $A$ be an arbitrary finite alphabet. The set of letters in a
string $u$ over $A$ is denoted by $L(u)$. String $u$ is called a
\emph{double occurrence string} if each $x \in L(u)$ occurs exactly
twice in $u$.  For example, $u = 41215425$ is a double occurrence
string over $L(u) = \{1,\ldots,5\}$. Let $\bar A = \{ \bar x \mid x
\in A \}$ with $A \cap \bar A = \emptyset$, and let $\tilde{A} = A
\cup \bar A$. We use the ``bar operator'' to move from $A$ to $\bar
A$ and back from $\bar A$ to $A$. Hence, for $x \in \tilde{A}$,
$\bar {\bar {x}} = x$. For a string $u = x_1 x_2 \cdots x_n$ with
$x_i \in A$, the \emph{inverse} of $u$ is the string $\bar u = \bar
x_n \bar x_{n-1} \cdots \bar x_1$.

We define the morphism $\pset{\cdot}:(\tilde{A})^* \rightarrow A^*$
as follows: for $x \in \tilde{A}$, $\pset{x} = x$ if $x \in A$, and
$\pset{x} = \bar{x}$ if $x \in \bar{A}$, i.e., $\pset{x}$ is the
``unbarred'' variant of $x$. Hence, e.g., $\pset{2 \bar 5 \bar 3} =
253$. A \emph{legal string} is a string $u \in (\tilde{A})^*$ where
$\pset{u}$ is a double occurrence string. We denote the empty string
by $\lambda$.

\begin{Example}
The string $u = q p s \bar q r p s r$ over $\tilde{A}$ with $A =
\{p,q,r,s\}$ is a legal string. As another example, the legal string
$3 4 4 5 6 7 5 6 7 8 9 \bar 3 \bar 2 2 8 9$ over $\tilde{B}$ with $B
= \{2,3,\ldots,9\}$ represents the micronuclear form of the gene
corresponding to the actin protein in the stichotrichous ciliate
\emph{Sterkiella nova}, see \cite{Biology_of_Ciliates,Ciliate_DB}.
\end{Example}

It is postulated that gene assembly is performed by three types of
elementary recombination operations, called loop, hairpin, and
double-loop recombination on DNA, see
\cite{Prescott_Ehren_Roz_2001}. These three recombination operations
have been modeled as three types of string rewriting rules operating
on legal strings \cite{Equiv_String_Graph_1,GeneAssemblyBook} ---
together they form the string pointer reduction system.  For all
$x,y \in \tilde{A}$ with $\pset{x} \not = \pset{y}$ we define:
\begin{itemize}
\item
the \emph{string negative rule} for $x$ by $\textbf{snr}_{x}(u_1 x x
u_2) = u_1 u_2$,
\item
the \emph{string positive rule} for $x$ by $\textbf{spr}_{x}(u_1 x
u_2 \bar x u_3) = u_1 \bar u_2 u_3$,
\item
the \emph{string double rule} for $x,y$ by $\textbf{sdr}_{x,y}(u_1 x
u_2 y u_3 x u_4 y u_5) = u_1 u_4 u_3 u_2 u_5$,
\end{itemize}
where $u_1,u_2,\ldots,u_5$ are arbitrary (possibly empty) strings
over $\tilde{A}$.

\begin{Example} \label{ex:string_contractions}
Let again $u = q p s \bar q r p s r$ be a legal string. We have
$\textbf{spr}_{q}(u) = \bar s \bar p r p s r$. And moreover,
$\textbf{spr}_{\bar r} \ \textbf{spr}_{\bar p} \ \textbf{spr}_{q}(u)
= \bar s \bar s$. Finally, $\textbf{snr}_{\bar s} \
\textbf{spr}_{\bar r} \ \textbf{spr}_{\bar p} \ \textbf{spr}_{q}(u)
= \lambda$.
\end{Example}

We now define a graph for a legal string representing whether or not
intervals within the legal string ``overlap''. Let $u = x_1 x_2
\cdots x_n$ be a legal string with $x_i \in \tilde{A}$ for $1 \leq i
\leq n$. For letter $y \in L(\pset{u})$ let $1 \leq i < j \leq n$ be
the positions of $y$ in $u$, i.e., $\pset{x_i} = \pset{x_j} = y$.
The \emph{$y$-interval} of $u$, denoted by $\mathrm{intv}_y$, is the
substring $x_k x_{k+1} \cdots x_l$ where $k = i$ if $x_i = y$ and $k
= i+1$ if $x_i = \bar y$, and similarly, $l = j$ if $x_j = y$ and $l
= j-1$ if $x_j = \bar y$ (i.e., a border of the interval is included
in case of $y$ and excluded in case of $\bar y$). Now the
\emph{overlap graph} of $u$, denoted by $\mathcal{G}_u$, is the
graph $(V,E)$ with $V = L(\pset{u})$ and $E = \{ \{x,y\} \mid x
\mbox{ occurs exactly once in } \pset{\mathrm{intv}_y}\}$. Note that
$E$ is well defined as $x$ occurring exactly once in the
$y$-interval of $u$ is equivalent to $y$ occurring exactly once in
the $x$-interval of $u$. Note that we have a loop $\{x\} \in E$ iff
both $x$ and $\bar x$ occur in $u$. The overlap graph as defined
here is an extension of the usual definition of overlap graph (also
called circle graph) from simple graphs (without loops) to graphs
(where loops are allowed). See
\cite[Section~7.4]{IntersectionGraphTheory_McKee_1999} for a brief
overview of (simple) overlap graphs.

\begin{Example}
The overlap graph $\mathcal{G}_u$ of $u = q p s \bar q r p s r$ is
exactly the graph $G$ of Example~\ref{ex:introd_pivot}.
\end{Example}

It is shown in \cite{Equiv_String_Graph_1,Equiv_String_Graph_2}, see
also \cite{GeneAssemblyBook}, that the string rules
$\textbf{snr}_{x}$, $\textbf{spr}_{x}$, and $\textbf{sdr}_{x,y}$ on
legal strings $u$ can be simulated as graph rules
$\textbf{gnr}_{x}$, $\textbf{gpr}_{x}$, and $\textbf{gdr}_{x,y}$ on
overlap graphs $\mathcal{G}_u$ in the sense that
$\mathcal{G}_{\textbf{spr}_{x}(u)} =
\textbf{gpr}_{x}(\mathcal{G}_u)$, where the left-hand side is
defined iff the right-hand side is defined, and similarly for
$\textbf{gdr}_{x,y}$ and $\textbf{gnr}_{x}$\footnote{There is an
exception for $\textbf{gnr}_{x}$: although
$\mathcal{G}_{\textbf{snr}_{x}(u)} =
\textbf{gnr}_{x}(\mathcal{G}_u)$ holds if the left-hand side is
defined, there are cases where the right-hand side is defined ($x$
is an isolated vertex in $\mathcal{G}_u$) while the left-hand side
is not defined ($u$ does not have substring $xx$). This is why the
string and graph models are ``almost'' equivalent. This difference
in models is not relevant for our purposes.}. It was shown in
\cite{BHH/PivotsDetPM/09} that ${\bf gpr}_{x}$ and ${\bf gdr}_{x,y}$
are exactly the two types of contractions of elementary pivots
$*\backslash\{x\}$ and $*\backslash\{x,y\}$ on a loop $\{x\}$ and an
edge $\{x,y\}$ without loops, respectively. The ${\bf gnr}_{x}$ rule
is the removal of isolated vertex $x$.

\begin{Example}
The sequence $\textbf{spr}_{\bar r} \ \textbf{spr}_{\bar p} \
\textbf{spr}_{q}$ applicable to $u$ given in
Example~\ref{ex:string_contractions} corresponds to a maximal
contraction of graph $G = \mathcal{G}_u$ of
Example~\ref{ex:introd_pivot} as can be seen in
Figure~\ref{fig:ex_max_contractions}.
\end{Example}

Within the theory of gene assembly one is interested in maximal
recombination strategies of a gene. These strategies correspond to
maximal contractions of a graph $G$ (hence decomposable into a
sequence $\varphi_1$ of contractions of elementary pivots ${\bf
gpr}$ and ${\bf gdr}$ applicable to (defined on) $G$) followed by a
sequence $\varphi_2$ of ${\bf gnr}$ rules, removing isolated
vertices, until the empty graph is obtained. Here we call these
sequences $\varphi = \varphi_2 \varphi_1$ of graph rules
\emph{complete contractions}. If we define the set of vertices $v$
of $\varphi$ used in ${\bf gnr}_{v}$ rules by $\gnrdom(\varphi)$,
then the following result holds by
Corollary~\ref{cor:max_contractions}.

\begin{Theorem} \label{thm:ga_gnr}
Let $G_1, G_2 \in \mathcal{O}_G$ for some graph $G$, and let
$\varphi$ be a complete contraction of graph $G_1$. Then there is a
complete contraction $\varphi'$ of $G_2$ such that $\gnrdom(\varphi)
= \gnrdom(\varphi')$.
\end{Theorem}

Hence, Theorem~\ref{thm:ga_gnr} shows that all the elements of
$\mathcal{O}_G$, for any graph $G$, have the same behaviour w.r.t.
the applicability of the rule ${\bf gnr}_{x}$.

A similar result as Theorem~\ref{thm:ga_gnr} was shown for the
string rewriting model, see
\cite[Theorem~34]{FibersRangeRedGraphs/Brijder07}\footnote{This
result states that two legal strings equivalent modulo ``dual''
string rules have the same reduction graph (up to isomorphism). It
then follows from \cite[Theorem~44]{StrategiesSnr/Brijder06} that
these strings have complete contractions with equal $\snrdom$, the
string equivalent of $\gnrdom$.}. It should be stressed however that
Theorem~\ref{thm:ga_gnr} is real generalization of the result in
\cite{FibersRangeRedGraphs/Brijder07} as not every graph has a
string representation (i.e., not every graph is an overlap graph),
and moreover it is obtained in a very different way: here the result
is obtained using techniques from linear algebra.

\section{Discussion}
We introduced the concept of dual pivot and have shown that it has
interesting properties: it has the same effect as the (regular)
pivot \emph{and} can be simulated by elementary row operations ---
consequently it keeps the kernel invariant. The dual pivot in this
way allows for an alternative definition of the (regular) pivot
operation. Furthermore, we have shown that two graphs have equal
kernel precisely when they have the same set of maximal pivots. From
this it follows that the set of maximal pivots is invariant under
dual pivot.

This main result is motivated by the theory of gene assembly in
ciliates in which maximal contractions correspond to complete
transformations of a gene to its macronuclear form. However, as
applying a maximal pivot corresponds to calculating a maximal
partial inverse of the matrix, the result is also interesting from a
purely theoretical point of view.

\subsection*{Acknowledgements}
We thank Lorenzo Traldi and the two anonymous referees for their
valuable comments on the paper. R.B. is supported by the Netherlands
Organization for Scientific Research (NWO), project ``Annotated
graph mining''.

\bibliography{../geneassembly}

\end{document}